\documentclass[letterpaper, 10 pt, conference]{ieeeconf}  % Comment this line out if you need a4paper

\IEEEoverridecommandlockouts                              % This command is only needed if
\usepackage{color}
\usepackage{graphics} % for pdf, bitmapped graphics files
\usepackage{epsfig} % for postscript graphics files
\usepackage{subfigure}
\usepackage{times} % assumes new font selection scheme installed
\usepackage{amsmath} % assumes amsmath package installed
\usepackage{amssymb}  % assumes amsmath package installed
\usepackage{lipsum}
\usepackage{multirow}

\newtheorem{theorem}{Theorem}[section]

\newtheorem{definition}{Definition}[section]
\newtheorem{lemma}[definition]{Lemma}

\newtheorem{corollary}[definition]{Corollary}
\newtheorem{remarkth}[definition]{Remark}
\newenvironment{remark}{\begin{remarkth}\upshape}{\end{remarkth}}

\usepackage{amssymb}

\newcommand{\proa}{A^*G \mbox{$\;$}_{\tau^*} \kern-3pt\times_\alpha
G \mbox{$\;$}_\beta \kern-3pt\times_{\tau^*} A^*G}

\title{\LARGE \bf
Variational obstacle avoidance problem on Riemannian manifolds
}

\author{Anthony Bloch$^{1}$, Margarida Camarinha$^{2}$ and Leonardo Colombo$^{1}$% <-this % stops a space
\thanks{$^{1}$ A. Bloch, and L. Colombo are with Department of Mathematics, University of Michigan, 530 Church St. Ann Arbor, 48109, Michigan, USA.
        {\tt\small abloch@umich.edu, ljcolomb@umich.edu}}%
\thanks{$^{2}$ M. Camarinha is with CMUC -- Centre for Mathematics of the University of Coimbra, Department of Mathematics, University of Coimbra
 3001-501 Coimbra, Portugal.
	    {\tt\small mmlsc@mat.uc.pt} }
}

\begin{document}

\maketitle
\thispagestyle{empty}
\pagestyle{empty}

%%%%%%%%%%%%%%%%%%%%%%%%%%%%%%%%%%%%%%%%%%%%%%%%%%%%%%%%%%%%%%%%%%%%%%%%%%%%%%%%
\begin{abstract}
We introduce variational obstacle avoidance problems on Riemannian manifolds and derive necessary conditions for the existence of their normal extremals. The problem consists of minimizing an energy functional depending on the velocity and covariant  acceleration, among a set of admissible curves, and also depending on a navigation function used to avoid an obstacle on the workspace, a Riemannian manifold.

We study two different scenarios, a general one on a Riemannian manifold and, a sub-Riemannian problem. By introducing a  left-invariant metric on a Lie group, we also study the variational obstacle avoidance problem on a Lie group. We apply the results to the obstacle avoidance problem of a planar rigid body and an unicycle.
\end{abstract}

%%%%%%%%%%%%%%%%%%%%%%%%%%%%%%%%%%%%%%%%%%%%%%%%%%%%%%%%%%%%%%%%%%%%%%%%%%%%%%%%
\section{Introduction}

Over the last few decades navigation functions have played a fundamental role in applications of trajectory planning for moving a system from a starting configuration to a goal configuration and, creating feasible and safe paths that avoid a prescribed obstacle minimizing some quantity such as energy or time.

In applications, navigation functions are typically given by artificial potential fields used for collision avoidance of certain regions through a radial analytic function on the configuration space \cite{Koditschek1990}. This approach has been studied by Khabit \cite{Khatib1986} for control problems and studied in the context of manifolds with boundary by Koditschek and Rimon \cite{Koditschek1990}. The mathematical foundations for the existence of smooth navigation functions on any smooth manifold have been proved by S. Smale \cite{Smale, HS}.% and M. Morse \cite{Morse}.

The theory of geodesics, presented, for instance, in Milnor \cite{Milnor}, is a  very rich example of the close relationship between variational problems and Riemannian geometry. Motivated by this connection and applications to dynamic interpolation on manifolds \cite{Jackson}, \cite{Noa:89}, Crouch and Silva Leite \cite{CroSil:91} started the development of an interesting geometric theory of generalized cubic  polynomials on a Riemannian manifold $M$, in particular on compact connected Lie groups endowed with a bi-invariant metric. Further extensions appear in the context of sub-Riemannian geometry, with connections with non-holonomic mechanics and control, studied by Bloch and Crouch \cite{BlCr, blochcrouch}. These sub-Riemannian problems are determined by additional constraints on a non-integrable distribution on $M$.

%\subsection{Main contributions:}

In this work we aim to introduce variational obstacle avoidance problems on $M$ as a first approach to further investigations related to dynamic interpolation and to avoid multiple regions in the workspace, as we explain at the end of the paper. We also aim to study necessary conditions for the existence of normal extremals in the variational problem among different situations, such as $M$ endowed with a left-invariant Riemannian metric on a Lie group, and sub-Riemannian problems where we must deal with constraints on a non-integrable distribution on $M$.

%\subsection{Outline:}
The structure of the paper is as follow. We start by introducing geometric structures on a Riemannian manifold that we will use together with admissible variation of curves and vector fields for the variational problem. Next, we introduce variational  obstacle avoidance problems on Riemannian manifolds and derive necessary conditions for the existence of normal extremals. In Section IV we extend our analysis to the sub-Riemannian situation where we also derive necessary conditions for the existence of normal extremals as in the general picture. By introducing a left-invariant Riemannian metric defined by an inner product on a Lie algebra of a Lie group we study the variational obstacle avoidance problem on a Lie group in Section V. We apply the results to obstacle avoidance problems for a planar rigid body and an unicycle in sections IV and V, respectively. Final comments and ongoing work are discussed at the end of the paper.
\section{Preliminaries on the calculus of variations}
\label{section2}
%In this section, we introduce the notion of a variational vector field. This
%allows us to introduce the required notation and tools to derive necessary
%conditions using the calculus of variations on manifolds (for a more detailed
%discussion, see, e.g.,\cite{MargaridaThesis}).

Let $M$ be a smooth ($\mathcal{C}^{\infty}$) \emph{Riemannian
manifold} with the \emph{Riemannian metric} denoted by
$\langle\cdot,\cdot\rangle:T_xM\times T_xM\to\mathbb{R}$ at each point $x\in M$, where $T_xM$ is
the \textit{tangent space} of $M$ at $x$. The length of a tangent vector is determined by its norm,
$||v_x||=\langle v_x,v_x\rangle^{1/2}$ with $v_x\in T_xM$.

A \emph{Riemannian connection} $\nabla$ on $M$, is a map that assigns to any two smooth vector fields $X$ and $Y$ on $M$ a new vector field, $\nabla_{X}Y$. For the properties of $\nabla$, we refer the reader to \cite{Boothby, bookBullo,Milnor}.  The operator
$\nabla_{X}$, which assigns to every vector field $Y$ the vector
field $\nabla_{X}Y$, is called the \emph{covariant derivative of
$Y$ with respect to $X$}. We denote by $\displaystyle{\frac{D}{dt}}$ the \textit{covariant time derivative}.

%Given a curve $x(t)$ and a vector field $X$, the covariant derivative of $X$ along $x$ is $\displaystyle{\frac{DX}{dt}=\nabla_{\dot{x}(t)}X}$. In local coordinates $x_1,\ldots,x_n$, the covariant derivative is $$(\nabla$$

Given vector fields $X$, $Y$
and $Z$ on $M$, the vector field $R(X,Y)Z$ given by \begin{equation}\label{eq:CurvatureTensorDefinition}
R(X,Y)Z=\nabla_{X}\nabla_{Y}Z-\nabla_{Y}\nabla_{X}Z-\nabla_{[X,Y]}Z
\end{equation}  is called the \emph{curvature tensor} of $M$. $[X,Y]$ denotes the \emph{Lie bracket} of the vector fields $X$ and $Y$. $R$ is trilinear in $X$, $Y$ and $Z$ and a tensor of type $(1,3)$. Hence for vector fields $X,Y,Z,W$ on $M$ the curvature tensor satisfies (\cite{Milnor}, p. 53) 
\begin{equation}\label{curvformula}\langle R(X,Y)Z,W\rangle=\langle R(W,Z)Y,X\rangle.\end{equation}

%The \emph{Lie bracket} of the vector fields $X$ and $Y$ denoted by $[X,Y]$ is defined by
%$[X,Y]f=X(Y f)-Y(X f)$, for all ${\cal C}^{\infty}$
%functions $f:M\to\mathbb{R}$.

\begin{lemma}[\cite{Boothby}, \cite{BlCr}]\label{lemma2}
Let $\omega$ be a one form on $(M,\langle\cdot,\cdot\rangle)$. The \textit{exterior derivative} of a one form $\omega$ is given by $$d\omega(X,Y)=X\omega(Y)-Y\omega(X)-\omega([X,Y])$$ for all vector fields $X,Y$ on $M$.

In particular, if $\omega(X)=\langle W,X\rangle$ it follows that
\begin{equation}\label{difomega}
d\omega(X,Y)=\langle\nabla_{X}W,Y\rangle-\langle\nabla_{Y}W,X\rangle.
\end{equation}
\end{lemma}

%\begin{lemma}[\cite{Boothby}]\label{lemma3} The curvature tensor $R$ on a Riemannian manifold $M$ satisfies the following properties for all vector fields $X,Y,Z,W,V$ on $M$.
%\begin{itemize}
%\item[(i)] $R(X,Y)Z+R(Y,X)Z=0$ and $\langle R(X,Y)Z,Z\rangle=0$.
%\item[(ii)] $R(X,Y)Z+R(Y,Z)X+R(Z,X)Y=0$ (1st Bianchi identity).
%\item[(iii)] $\langle R(X,Y)Z,W\rangle=\langle R(W,Z)Y,X\rangle$.
%\item[(iv)] As a consequence of (i) and (iii), $\langle R(X,Y)Z,W\rangle+\langle R(X,Y)W,Z\rangle=0.$
%\item[(v)] $\langle(\nabla_YR)(X,Z)W,V\rangle=\langle (\nabla_YR)(V,W)Z,X\rangle$.
%\end{itemize}
%\end{lemma}

Let $\Omega$ be the set of all ${\cal C}^1$ piecewise smooth curves $x:[0, T]\to M$ in
$M$ such that $x(0)$, $x(T)$, $\frac{dx}{dt}(0)$ and $\frac{dx}{dt}(T)$ are fixed. The set $\Omega$ is called the
\textit{admissible set}.

For the class of curves in $\Omega$, we introduce the ${\cal C}^1$ piecewise smooth \textit{one-parameter admissible variation} of a curve $x\in\Omega$ by $
\alpha : (-\epsilon , \epsilon ) \times [0,T]  \rightarrow  M ;(r,t) \mapsto \alpha (r,t)=\alpha_r(t)$
that verify $\alpha_0=x$ and $\alpha_r\in\Omega$, for each $r\in (-\epsilon , \epsilon )$.

The \textit{variational vector field}  associated to an admissible variation $\alpha$ is a ${\cal C}^{1}$-piecewise smooth vector fields along $x$ defined by $$X(t)=\frac{D}{\partial r}\Big{|}_{r=0}\alpha(r,t)\in T_{x}\Omega,$$
verifying  the boundary conditions
\begin{equation}\label{3.6}
X(0)=0,\quad   X(T)=0,\quad
\frac{DX}{dt}(0)=0,\quad   \frac{DX}{dt}(T)=0,\end{equation} where the tangent space  of  $\Omega$ at $x$ is the
 vector space $T_{x}\Omega$  of all ${\cal C}^{1}$ piecewise smooth vector
fields $X$ along $x$  verifying the boundary conditions (\ref{3.6}).

Consider a vector field $W$  along a curve $x$ on $M$. The $s$th-order covariant  derivative along $x$  of $W$ is denoted by $\displaystyle{\frac{D^{s}W}{dt^{s}}}$, $s\geq 1$. We also denote by $\displaystyle{\frac{D^{s+1}x}{dt^{s+1}}}$ the $s$th-order covariant derivative along $x$ of the velocity vector field of $x$, $s\geq 1$.

\begin{lemma}[\cite{Milnor}, p.$52$]\label{lemma} The one-parameter variation satisfies
$$\frac{D}{\partial r}\frac{D^2\alpha}{\partial
t^2}=\frac{D^2}{dt^2}\frac{\partial \alpha}{\partial r}+R\left(\frac{\partial \alpha}{\partial r}
,\frac{\partial \alpha}{\partial t}\right)\frac{\partial
\alpha}{\partial t}$$ where $R$ is the curvature tensor.
\end{lemma}
%{\color{red} This lemma is not correct like that and it is not necessary to use it. I suggest to take it off.
%\begin{lemma}[see, for instance, \cite{MargaridaThesis}]\label{lemma4} The following identity for the differentiation of the curvature tensor holds
%\begin{equation*}
%\frac{D}{dt}\left(R\left(\frac{\partial\alpha}{\partial r},\frac{\partial\alpha}{\partial t}\right)\frac{\partial\alpha}{\partial t}\right)=-R\left(\frac{D^{2}\alpha}{\partial t^2},\frac{\partial\alpha}{\partial t}\right)\frac{\partial\alpha}{\partial r}.
%\end{equation*}
%\end{lemma}
%}
%
\section{The variational obstacle avoidance problem on a Riemannian manifold}

%Let $(M,\langle\cdot,\cdot\rangle)$ be an $n$-dimensional Riemannian manifold where $\langle\cdot,\cdot\rangle$ denotes the Riemannian metric. The associated Levi-Civita connection on $M$ is denoted as $\nabla $ and $R$ denotes the curvature tensor of  $\nabla$.

Let $T$, $\sigma$ and $\tau$  be positive real numbers, $(p_0,v_0)$, $(p_T,v_T)$ points in $TM$  and $q$ a point on $M$ representing an obstacle on the workspace $M$. Consider the set $\Omega$ of all ${\cal C}^{1}$ piecewise smooth curves  on $M$, $x:[0,T]\rightarrow M$ verifying the boundary conditions \begin{equation}\label{3.1}
x(0)=p_0, \quad  x(T)=p_T,\quad \frac{dx}{dt}(0)=v_0, \quad \frac{dx}{dt}(T)=v_T,\end{equation} and define the functional $J$ on $\Omega$ given by
\begin{equation}\label{3.2}
J(x)=\int_{0}^{T}\frac{1}{2}\left(\Big{\|}\frac{D^2x}{dt^2}(t)\Big{\|}^2+
 \sigma \Big{\|}\frac{ dx}{dt}(t)\Big{\|}^2+V_q(x(t))\right)dt.
\end{equation}

\noindent This functional is constructed as a linear combination of the velocity and the covariant acceleration of the trajectory regulated by the
parameter $\sigma$, together with a navigation function used to avoid the obstacle $q$ described as the zero level surface of a know scalar valued analytic function (see, e.g., \cite{Khatib1986}, \cite{K88}, \cite{Koditschek1990}).

The navigation function $V_q$ is an artificial potential field-based function
represented by a force inducing an artificial repulsion from the surface of the obstacle. We use the approach introduced by Khatib \cite{Khatib1986} which consists on using a local inverse potential field going to infinity as the inverse square of a know scalar valued analytic function near the obstacle, and decay to zero at some positive level surface far away from the obstacle, in order that a particle on such a field never hits $q$.

%, created by a potential field.

%nverse square distance between the obstacle $q\in M$ and the trajectory position $x\in M$ . It is given by $$V_q(x)=\frac{\tau}{d^2(q,x)},$$ where the parameter $\tau$ allows to control the strenght of the  potential f\textcolor{blue}{field}.

\textbf{Problem}: The \textit{variational obstacle avoidance problem} consists in minimizing the functional $J$ among $\Omega$.

In order to minimize the functional $J$ among  the set $\Omega$ we want to find curves $x\in \Omega$ verifying $J(x)\leq J(\tilde{x})$, for all admissible curves $\tilde{x}$ in a ${\cal C}^1$-neighborhood of $x$.
%We call this problem the {\em obstacle avoidance variational problem}.

To compare the value of $J$ at a curve $x\in \Omega$ to the value of $J$ at a
nearby curve $\tilde{x}\in \Omega$, we use one-parameter admissible
variations $\alpha$ of $x\in \Omega$.%, that is,
%one-parameter variations
%$
%\alpha : (-\epsilon , \epsilon ) \times [0,T]  \rightarrow  M ;(r,t) \mapsto \alpha (r,t)=\alpha_r(t)$
%that verify $\alpha_0=x$ and $\alpha_r\in\Omega$, for each $r\in (-\epsilon , \epsilon )$.
%The variational vector field  associated to an admissible variation $\alpha$, $X=\frac{\partial \alpha}{ \partial r}_{|r=0}$,
%is a ${\cal C}^{1}$ piecewise smooth vector
%fields along $x$
%verifying  the boundary conditions
%\begin{equation}\label{3.6}
%\begin{array}{ll}
% \displaystyle X(0)=0,&  \displaystyle  X(T)=0,\\
%\\
% \displaystyle \frac{DX}{dt}(0)=0,&  \displaystyle   \frac{DX}{dt}(T)=0.\\
%\\
%\end{array}
%\end{equation}
%
%\noindent  The tangent space  of  $\Omega$ at $x$ is the
% vector space $T_{x}\Omega$  of all ${\cal C}^{1}$ piecewise smooth vector
%fields $X$ along $c$  verifying the boundary conditions (\ref{3.6}).
%
%\par
\begin{theorem} \label{t3.2}
	Let $x\in \Omega$.
 If $\alpha$ is an admissible variation of $x$ with variational vector field
$X \in T_x\Omega$, then \begin{align*}\label{3.7}
&\frac{d}{dr}J(\alpha _{r})\Big{|}_{r=0}=\int_{0}^{T}\Big{\langle}X,\frac{D^{4}x}{dt^{4}}+R\left(\frac{D^2x}{dt^2},\frac{dx}{dt}\right)\frac{dx}{dt}\\
&\qquad\qquad\qquad\qquad\qquad\qquad-\sigma \frac{D^2x}{dt^2}+\frac{1}{2}\mbox{grad } V_q(x)\Big{\rangle} dt \\
&\qquad\qquad\qquad+\sum_{i=1}^{l-1} \Big{\langle}\frac{DX(t_{i})}{dt},\frac{D^2x}{dt^2}(t^{+}_{i})-
\frac{D^2x}{dt^2}(t^{-}_{i})\Big{\rangle}\\
&\qquad\qquad\qquad- \sum_{i=1}^{l-1} \Big{\langle}X(t_{i}),\frac{D^{3}x}{dt^{3}}(t^{+}_{i})-
\frac{D^{3}x}{dt^{3}}(t^{-}_{i})\Big{\rangle}.
\end{align*}
%\end{equation}
\end{theorem}

\textbf{Proof}: Let $\alpha$ be an admissible variation of $x$ with variational vector field
$X \in T_x\Omega$. Then
\begin{align*}
\frac{d}{dr}J(\alpha _{r})=&\int_{0}^{T}\left(\Big{\langle}\frac{D}{\partial
r}\frac{D^2 \alpha }{\partial t^2},\frac{D^2 \alpha }{\partial t^2}\Big{\rangle}+\sigma\Big{\langle}\frac{D^2\alpha}{\partial r \partial t},\frac{\partial \alpha }{\partial
t}\Big{\rangle}\right.\\
&\left.\qquad\qquad\qquad+ \frac{1}{2}\frac{\partial}{\partial r}V(q,\alpha)\right)dt.
\end{align*}

By considering the gradient vector field $(\mbox{grad }V_q)$ of the potential field-based function $V_q: M\to \mathbb{R}$ we have $$ \frac{\partial}{\partial
r}V_q(\alpha)=\Big{\langle}\frac{\partial \alpha}{\partial r },\mbox{grad }V_q(\alpha)\Big{\rangle}.$$

%\noindent By other hand, we know that
%$$
%\begin{array}{c}
%\displaystyle \frac{D}{\partial r}\frac{D^2\alpha}{\partial
%t^2}=\frac{D^2}{dt^2}\frac{\partial \alpha}{\partial r}+R(\frac{\partial \alpha}{\partial r}
%,\frac{\partial \alpha}{\partial t})\frac{\partial
%\alpha}{\partial t}.\\
%\\
% \end{array}
%$$

%\noindent The two identities stated above imply that
By Lemma \ref{lemma} and the previous identity we have
\begin{align*}
\frac{d}{dr}J(\alpha _{r})=&\int_{0}^{T}\left(
\Big{\langle}\frac{D^2}{dt^2}\frac{\partial \alpha}{\partial r},\frac{D^2 \alpha }{\partial
t^2}\Big{\rangle}\right.\\
&\qquad+\Big{\langle}R\left(\frac{\partial \alpha}{\partial r},\frac{\partial
\alpha}{\partial t}\right)\frac{\partial \alpha}{\partial t}, \frac{D^2 \alpha }{\partial t^2}\Big{\rangle}\\
&\left.\qquad+\sigma
\Big{\langle}\frac{D^2\alpha}{\partial t \partial r},\frac{\partial \alpha }{\partial
t}\Big{\rangle}+\Big{\langle}\frac{\partial \alpha}{\partial r }\frac{1}{2}\mbox{grad }V_q(\alpha)\Big{\rangle}\right)dt.
\end{align*}

Integrating the first term by parts twice, the third term once,
and applying the property (\ref{curvformula}) of the curvature tensor $R$ to the second term, we obtain
\begin{align*}
\frac{d}{dr}J(\alpha _{r})= &
\sum_{i=1}^{l}\left[\Big{\langle}\frac{D}{\partial t}\frac{\partial \alpha}{\partial r},\frac{D^2 \alpha
}{\partial t^2}\Big{\rangle}-\Big{\langle}\frac{\partial \alpha}{\partial r}, \frac{D^3 \alpha }{\partial t^3}\Big{\rangle}\right.\\
&\left.+\sigma\Big{\langle}\frac{\partial \alpha }{\partial r},\frac{\partial \alpha }{\partial
t}\Big{\rangle}\right]_{t_{i-1}^+}^{t_i^-}\\
&+\int_{0}^{T}\left(
\Big{\langle}\frac{\partial \alpha}{\partial r},\frac{D^4 \alpha }{\partial
t^4}+R\left(\frac{D^2 \alpha}{\partial t^2},\frac{\partial
\alpha}{\partial t}\right)\frac{\partial \alpha}{\partial t}\right.\\
&\left.-\sigma\frac{D^2 \alpha }{\partial t^2}+\frac{1}{2}\mbox{grad }V_q(\alpha)\Big{\rangle}\right)dt,
\end{align*} where given that $x\in\Omega$ we must consider a partition of the interval $[0,T]$ as $0=t_0<t_1<\ldots<t_l=T$ in a way that
$x$ is smooth in each subinterval. % divide the integral in each subinterval determined by the partition. %with $\displaystyle{t_i^{-}=\lim_{t\to t_i^{-}}t}$ and $\displaystyle{t_i^{+}=\lim_{t\to t_i^{+}}t}$.

Next, by taking $r=0$ in the last equality, we obtain
\begin{align*}
\frac{d}{dr}J(\alpha _{r})\Big{|}_{r=0}=&\sum_{i=1}^{l} \left[\Big{\langle}\frac{DX}{dt},\frac{D^2x}{dt^2}\Big{\rangle}\right.\\
&\qquad\left.- \Big{\langle}X,\frac{D^{3}x}{dt^{3}}\Big{\rangle}+\sigma\Big{\langle} X,\frac{dx}{dt}\Big{\rangle}\right]_{t_{i-1}^{+}}^{t_i^{-}}\\
+&\int_{0}^{T}\left(\Big{\langle}X,\frac{D^{4}x}{dt^{4}}+R\left(\frac{D^2x}{dt^2},\frac{dx}{dt}\right)\frac{dx}{dt}\right.\\
&\left.\quad\qquad-\sigma \frac{D^2x}{dt^2}+\frac{1}{2}\mbox{grad } V_q(x(t))\Big{\rangle}\right) dt .
\end{align*}

 %\displaystyle \frac{d}{dr}J(\alpha _{r})_{|r=0} & \displaystyle = &\displaystyle
%\sum_{i=1}^{l}\left[ <\frac{DX}{d t},\frac{D^2x
%}{d t^2}>-<X, \frac{D^3 x}{d t^3}>+\sigma<X,\frac{dx}{d
%t}>\right]_{t_{i-1}^+}^{t_i^-}\\
%&&\\
%&& \displaystyle
%+\int_{0}^{T}(
%<X,\frac{D^4 x}{d
%t^4}+R(\frac{D^2 x}{d t^2},\frac{d
%x}{d t})\frac{dx}{d t}-\sigma\frac{D^2x }{d t^2}+ \mbox{grad }V_q(x)>)dt. \\
%&&
%\\
 %\end{array}
%$$

Since the vector field $X$ is ${\cal C}^1$, piecewise smooth
on $[0,T]$, verifies the boundary conditions (\ref{3.6}) and the curve $x$
is ${\cal C}^1$ on $[0,T]$, the result follows.
\quad$\Box$

\vspace*{.3cm}

\par
 \begin{theorem}\label{t3.3}
If  $x \in \Omega$ is a local minimizer of
$J$, then $x$ is smooth on $[0,T]$ and verifies
\begin{equation} \label{3.8}
\frac{D^{4}x}{dt^{4}}+R\left(\frac{D^2x}{dt^2},\frac{dx}{dt}\right)\frac{dx}{dt}- \sigma
\frac{D^2x}{dt^2}+ \frac{1}{2}\mbox{grad }V_q(x) \equiv 0.
\end{equation}
\end{theorem}
\vspace{.5cm}

\textbf{Proof.} Assume $x \in \Omega$ is a local minimizer of $J$
over $\Omega$. Then  $\displaystyle{\frac{d}{dr}J(\alpha_r)\mid_{r=0}=0}$, for each admissible variation $\alpha$  of $x$ with variational vector field
$X \in T_x\Omega$.

Let us consider  $X \in T_x\Omega$ defined by
$$f\left[\frac{D^{4}x}{dt^{4}}+R\left(\frac{D^2x}{dt^2},\frac{dx}{dt}\right)\frac{dx}{dt}- \sigma
\frac{D^2x}{dt^2}+\frac{1}{2}\mbox{grad }V_q(x)\right],
$$ where $f$ is a smooth real-valued function on $[0,T]$ verifying
$f(t_i)=f^{\prime}(t_i)=0$ and $f(t)>0$, $t\neq t_i$, $i=1,\ldots, l-1$.
So, we have
\begin{align*}
\displaystyle \frac{d}{dr}J(\alpha_r)\Big{|}_{r=0} =&\int_0^Tf(t)\Big{|}\Big{|}\frac{D^{4}x}{dt^{4}}+R\left(\frac{D^2x}{dt^2},\frac{dx}{dt}\right)\frac{dx}{dt}\\
&\quad\qquad\qquad- \sigma
\frac{D^2x}{dt^2}+\frac{1}{2}\mbox{grad }V_q(x) \Big{|}\Big{|}^2dt
\end{align*} and since $f(t)>0$ for $t\in [0,T]$, except in a set of measure zero, it follows that
$$\Big{|}\Big{|} \frac{D^{4}x}{dt^{4}}+R\left(\frac{D^2x}{dt^2},\frac{dx}{dt}\right)\frac{dx}{dt}- \sigma
\frac{D^2x}{dt^2}+\frac{1}{2}\mbox{grad }V_q(x)\Big{|}\Big{|}\equiv 0,$$

\noindent which leads to the equation (\ref{3.8}).

Next let us choose the vector field $X \in T_x\Omega$ so that
$$X(t_{i})=\frac{D^{3}x}{dt^{3}}(t^{+}_{i})-
\frac{D^{3}x}{dt^{3}}(t^{-}_{i})$$

\noindent and
$$\frac{DX(t_{i})}{dt}=\frac{D^2x}{dt^2}(t^{-}_{i})-
\frac{D^2x}{dt^2}(t^{+}_{i}),$$

\noindent for $i=1,\ldots, l-1$.
Thus,
\begin{align*}\frac{d}{dr}J(\alpha_r)\Big{|}_{r=0}=&\sum_{i=1}^{l-1}\left(\Big{|}\Big{|}\frac{D^2x}{dt^2}(t^{+}_{i})-
\frac{D^2x}{dt^2}(t^{-}_{i})\Big{|}\Big{|}^2 \right.\\
&\left.+ \Big{|}\Big{|} \frac{D^{3}x}{dt^{3}}(t^{+}_{i})-
\frac{D^{3}x}{dt^{3}}(t^{-}_{i})\Big{|}\Big{|}^2\right)=0,\end{align*}

\noindent which implies that
$$\displaystyle  \frac{D^2x}{dt^2}(t^{+}_{i})=
\frac{D^2x}{dt^2}(t^{-}_{i})  \; \; {\rm and}\; \; \frac{D^{3}x}{dt^{3}}(t^{+}_{i})=
\frac{D^{3}x}{dt^{3}}(t^{-}_{i}) .$$

\noindent Hence, $x$ is is smooth
on $[0,T]$.
\quad$\Box$
\begin{remark}
When $V=0$, equation (\ref{3.8}) reduces to the cubic polynomials in tension equation \cite{SCC00}
\begin{equation}\label{eq11}
\frac{D^{4}x}{dt^{4}}+R\left(\frac{D^2x}{dt^2},\frac{dx}{dt}\right)\frac{dx}{dt}- \sigma
\frac{D^2x}{dt^2} \equiv 0.
\end{equation}

\end{remark}

\section{Sub-Riemannian variational problem}

Next, we extend our analysis to the sub-Riemannian context, that is, we assume the velocity vector field $\displaystyle{\frac{dx}{dt}}$ lies on some
distribution $\mathcal{D}\subset TM$. This distribution $\mathcal{D}$ is defined by
non-integrable constraints on the velocity vector field determined by one-forms
$\omega _j\in T^{*}M$ with $1\leq j \leq k<n$, satisfying
\begin{equation}\label{eqconstraint}
\omega_j\left(\frac{dx}{dt}\right)=\Big{\langle}Y_j,\frac{dx}{dt}\Big{\rangle}=0,
\end{equation} where $Y_1, \cdots ,Y_k,\cdots , Y_n$ are linearly independent vector fields on some neighborhood $\Omega$ of $x\in M$.% and where $c_j$, $ 1\leq j\leq k<n$ are given constants.

To deal with the constraints we also need to define the tensors $S_{i}$, $(S_{i})_{x}:T_{x}M\to T_{x}M$ by $$d\omega_j(u,z)=<S_j(u),z>=-<S_j(z),u>, u,z \in T_xM$$
%$$d\omega_{i}(X,Y)=\langle S_{i}(X),Y\rangle=-\langle S_{i}(Y),X\rangle$$ where $i=1,\ldots,k$ and $d$ is the exterior derivative.

\textbf{Problem}: The sub-Riemannian variational obstacle avoidance problem consists in minimizing the functional $J$ defined on \eqref{3.1} among $\Omega$ with the additional constraints \eqref{eqconstraint}.

This type of problem was studied in Bloch and Crouch \cite{BlCr} and Crouch and Silva Leite
\cite{CroSil:95}.

 We derive necessary conditions for the existence of normal extremals in this sub-Riemannian problem, by extending our previous analysis for the general case following the result of Bloch and Crouch \cite{BlCr}, \cite{blochcrouch}.
\begin{theorem} \label{T4}
A necessary condition for $x\in\Omega $ to be a normal extremal for the sub-Riemannian variational obstacle avoidance problem is that $x$ be of class ${\cal C}^2$ and there exist smooth
functions $\lambda _j$, $j=1,\cdots ,k$ (the Lagrange multipliers) such
that, for every $t\in [t_{i-1},t_i],\; i=1,\cdots ,l$, the following equations holds
\begin{align*}\label{eq:23}
0=&\frac{D^{4}x}{dt^{4}}+R\left(\frac{D^2x}{dt^2},\frac{dx}{dt}\right)\frac{dx}{dt}- \sigma
\frac{D^2x}{dt^2}+\frac{1}{2}\mbox{grad }V_q(x)\\
&-\sum_{j=1}^k\lambda _j^{\prime}Y_j-\sum_{j=1}^k\lambda _jS_j\left(\frac{dx}{dt}\right),
\end{align*} together with $ \Big{\langle}Y_j,\frac{dx}{dt}\Big{\rangle}=0, \; 1\leq j\leq k$.

\end{theorem}

\textbf{Proof}: Consider the extended functional
\begin{align*}
\widetilde{J}(x)=&\frac{1}{2}\int_{0}^{T}\left(\Big{\|}\frac{D^2x}{dt^2}(t)\Big{\|}^2+
 \sigma \Big{\|}\frac{dx}{dt}(t)\Big{\|}^2+V_q(x(t))\right.\\
 &\left.\qquad\qquad+\sum_{j=1}^{k}\lambda_j\Big{\langle}Y_j,\frac{dx}{dt}\Big{\rangle}\right)dt.
\end{align*}

We derive necessary conditions for existence of normal extremals by studying the equation $$\frac{d}{dr}\tilde{J}(\alpha_r)\Big{|}_{r=0}=0$$ for $\alpha$ an admissible variation of $x$ with variational vector field $X \in T_x\Omega$ and $\lambda_j$ the Lagrange multipliers.

Taking into account the proof of Theorems \ref{t3.2} and \ref{t3.3} we only need to study the influence of variations in the term $\displaystyle{\sum_{j=1}^{k}\lambda_j\Big{\langle}Y_j,\frac{dx}{dt}\Big{\rangle}}$ where the vector fields $Y_{j}$ on $M$ are determined by $\omega_{j}(Z)=\langle Y_{j}, Z\rangle$, $j=1,\ldots,k$ for each vector field  $Z$ on $M$. Therefore, $\displaystyle{\frac{d}{dr}\tilde{J}(\alpha_r)\Big{|}_{r=0}}$ must have two additional terms compared with $\displaystyle{\frac{d}{dr}J(\alpha_r)\Big{|}_{r=0}}$. Those terms are $$\sum_{j=1}^{k}\lambda_{j}\Big{\langle}\nabla_{\frac{\partial\alpha}{\partial r}}Y_{j},\frac{\partial\alpha}{\partial t}\Big{\rangle}+\sum_{j=1}^{k}\lambda_{j}\Big{\langle}Y_{j},\frac{D^2\alpha}{\partial t \partial r}\Big{\rangle}.$$

After integration by parts in the second term and evaluating at $r=0$, the integrand can be re-written with the additional terms $$\sum_{j=1}^{k}\lambda_j\langle\nabla_{X}Y_{j},\frac{dx}{dt}\rangle-\sum_{j=1}^{k}\lambda_{j}'\langle Y_{j},X\rangle-\lambda_j\sum_{j=1}^{k}\langle\frac{DY_{j}}{dt},X\rangle.$$ Using the identity (\ref{difomega})  the new terms compared with the ones provided by Theorems \ref{t3.2} and \ref{t3.3} which give rise to necessary conditions for the existence of normal extremals in this sub-Riemannian problem are: $$-\sum_{j=1}^{k}\lambda_{j}d\omega_j\left(\frac{dx}{dt},X\right)-\sum_{j=1}^{k}\lambda_j'\Big{\langle}Y_j,X\Big{\rangle}.$$ Using the fact that $d\omega_j\left(\frac{dx}{dt},X\right)=\Big{\langle}S_j\left(\frac{dx}{dt}\right),X\Big{\rangle}$ the result follows. \quad$\Box$
%\section{The special euclidean group $SE(2)$}

\begin{corollary}
Any abnormal extremal for the sub-Riemannian variational obstacle avoidance problem satisfy $$\sum_{j=1}^k\lambda _j^{\prime}Y_j+\sum_{j=1}^k\lambda _jS_j\left(\frac{dx}{dt}\right)=0$$ where $\lambda_j$, $j=1,\ldots,k$ are not all identically zero.\end{corollary}

\subsection{Application to variational obstacle avoidance problem for a planar rigid body on $SE(2)$.}
The  special euclidean Lie group $SE(2)$ consists of all the transformations of $\mathbb{R}^2$ of the form $z \mapsto Rz+v$, where $v\in \mathbb{R}^2$ and $R\in SO(2)$. This Lie group is
isomorphic to the semidirect  product Lie group $SO(2)\ltimes \mathbb{R}^2$. The transformations can be represented by $(R, v)$, where
$$R=\left(
     \begin{array}{cc}
       \cos \theta & -\sin \theta\\
       \sin \theta & \cos \theta \\
     \end{array}
   \right)$$
or, for the sake of simplicity, by the matrix
$$\left(
     \begin{array}{cc}
      R & v\\
       0 & 1 \\
     \end{array}
   \right)=\left(
     \begin{array}{ccc}
       \cos \theta & -\sin \theta & x\\
       \sin \theta & \cos \theta & y\\
       0&0& 1
     \end{array}
   \right).$$
The composition law is defined by
$(R,v)\cdot(S,w)=(RS,Rw+v)$
with identity element $(I,0)$ and inverse
$(R,v)^{-1}=(R^{-1},-R^{-1}v)$.

%Note that the composition law corresponds to the usual matrix multiplication if we consider the matrix representation.
The Lie algebra $\mathfrak{se}(2)$ of $SE(2)$ is determined by
$$\mathfrak{se}(2)=\Big{\{}\left(
     \begin{array}{cc}
      A & b\\
       0 & 0 \\
     \end{array}
   \right): A\in \mathfrak{so}(2) \hbox{ and } b\in \mathbb{R}^2\Big{\}}.$$
For simplicity, we write $A=-aJ$ where $J=\left(
     \begin{array}{cc}
      0& 1\\
       -1 &0 \\
     \end{array}
   \right)$ and we identify the Lie algebra $\mathfrak{se}(2)$ with $\mathbb{R}^3$ via the isomorphism
   $\displaystyle{\left(
     \begin{array}{cc}
     -aJ& b\\
       0 & 0 \\
     \end{array}
   \right)\mapsto (a,b)}$.

   The Lie bracket in $\mathbb{R}^{3}$ is given by
$[(a,b),(c,d)]=(0,-aJd+cJb).$
The  basis of  $\mathfrak{se}(2)$ represented by the canonical basis of $\mathbb{R}^3$ verifies
$[e_1,e_2]=e_3$, $[e_2,e_3]=0$, $[e_3,e_1]=e_2.$
%The adjoint action is given by
%$$Ad_{(R,v)}(a,b)=(a,aJv+Rb).$$

%%%%%%%%%%%%%%%%%%%

%\textbf{\large Reformulei a vers\~{a}o anterior a partir daqui. Comecei por considerar o caso mais simples do corpo r\'{\i}gido planar.}

%%%%%%%%%%%%%%%%%%%%%%%%%%

The Riemannian metric on $SE(2)\simeq\mathbb{R}^{2}\times S^{1}$, locally parametrized by $\gamma=(x,y,\theta)$, is determined by the matrix $\mbox{diag }(m,m,J)$. %($J=2$ and $m=1$)
%The Christoffel symbols of the Levi-Civita connection are zero.
 The curvature tensor is zero.
%and the metric is specified with $J=2$, $m=1$
We consider the navigation function \begin{equation}\label{navfunct}V(\gamma)=\frac{\tau}{x^2+y^2-1}\end{equation} representing an obstacle with circular shape and unitary radius in the $xy$-plane, centered at the origin, with $\tau\in\mathbb{R}^{+}$. %By considering $$\frac{dx}{dt}=\theta^{\prime}\frac{\partial}{\partial \theta}+x^{\prime}\frac{\partial}{\partial x}+y^{\prime}\frac{\partial}{\partial y},$$ one obtain $$\frac{D^{i}x}{dt^{i}}=\theta^{(i)}\frac{\partial}{\partial \theta}+x^{(i)}\frac{\partial}{\partial x}+x^{(i)}\frac{\partial}{\partial y}.$$
%Moreover,
Note that $$\mbox{grad } V_q(\gamma)=-\frac {2\tau}{m(x^2+y^2-1)^2}\left(x\frac{\partial}{\partial x}+y\frac{\partial}{\partial y}\right).$$
%Therefore, equations

%$$\frac{D^{4}x}{dt^{4}}+R\left(\frac{D^2x}{dt^2},\frac{dx}{dt}\right)\frac{dx}{dt}- \sigma
%\frac{D^2x}{dt^2}+\mbox{grad }V_q(x)=0$$

%becomes in

%\begin{align*}
%0=&(\theta^{(4)}-\sigma\theta^{(2)})\frac{\partial}{\partial \theta}+(x^{(4)}-\sigma x^{(2)}-\frac {\tau x}{(x^2+y^2)^2})\frac{\partial}{\partial x}\\&+(y^{(4)}-\sigma y^{(2)}-\frac {\tau y}{(x^2+y^2)^2})\frac{\partial}{\partial y}.\end{align*}

By Theorem \ref{t3.3}, the equations determining necessary conditions for the existence of normal extremals in the variational problem are\begin{align*}
                   \theta^{(4)}&=\sigma\theta'',\\
                  x^{(4)}&=\sigma x''+\frac {\tau x}{m(x^2+y^2-1)^2},\\
                  y^{(4)}&=\sigma y''+\frac {\tau y}{m(x^2+y^2-1)^2},
                 \end{align*}  with given boundary conditions

             $(x(0),y(0),\theta(0)),\qquad (x(T),y(T),\theta(T))$

                 $(x^{\prime}(0),y^{\prime}(0),\theta^{\prime}(0)), \quad (x^{\prime}(T),y^{\prime}(T),\theta^{\prime}(T))$.

\subsection{Application to obstacle avoidance sub-Riemannian problem for an unicycle.}

We study motion planning of a unicycle with obstacles. To avoid the obstacle, we use the navigation function approach.

The unicycle is a homogeneous disk on a horizontal plane and it is equivalent to a wheel rolling on a plane \cite{Bl, bookBullo}. The configuration of the unicycle at any given time is completely determined by the element $(R,v)\in\mathrm{SE}(2)\cong\mathbb{R}^{2}\times\mathrm{S}^{1}$. As before, we consider $SE(2)$ locally parametrized by $\gamma=(x,y,\theta)$ and also consider the navigation function $V(\gamma)$ given in \eqref{navfunct} representing a circular obstacle of unitary radius in the $xy$-plane centered at the origin.

The distribution is spanned by the one-form $\omega=\sin \theta dx-\cos \theta dy$ with corresponding vector field given by (see \cite{bookBullo} for instance) $$Y_1= \frac 1 m \left(\sin \theta \frac{\partial}{\partial x}-\cos \theta \frac{\partial}{\partial y}\right).$$ Note that
$$S(U)=-\frac 1J(u_2\cos \theta +u_3\sin \theta)\frac{\partial}{\partial \theta}+\frac {u_1}m(\cos \theta  \frac{\partial}{\partial x}+\sin \theta \frac{\partial}{\partial y}),$$
for each  vector field on $SE(2)$ denoted by $U=u_1 \frac{\partial}{\partial \theta}+ u_2 \frac{\partial}{\partial x}+u_3 \frac{\partial}{\partial y}$.
%$$S(u)=-\frac 1J(\cos \theta +\sin \theta)\frac{\partial}{\partial \theta}+\frac {1}m(\cos \theta  \frac{\partial}{\partial x}+\sin \theta \frac{\partial}{\partial y}).$$
 We complete $\hbox{span}\{Y_1\}$ to a basis of vector fields by consider
 $\displaystyle{Y_2=\frac{1}{J}\frac{\partial}{\partial \theta}}$, $\displaystyle{Y_3=\frac{\cos \theta}{m} \frac{\partial}{\partial x}+\frac{\sin \theta}{m}\frac{\partial}{\partial y}}.$

By Theorem \ref{T4} the equations determining necessary conditions for the existence of normal extremals in the variational problem are
\begin{align*}
                    \theta^{(4)}&=\sigma\theta'' -\frac 1J \lambda(x^{\prime}  \cos  \theta+ y^{\prime} \sin  \theta), \\
                  x^{(4)}&=\sigma x''+\frac {\tau x}{m(x^2+y^2-1)^2} +\frac 1m\lambda^{\prime}\sin \theta+\frac 1m\lambda\theta^{\prime} \cos \theta,\\
                  y^{(4)}&=\sigma y''+\frac {\tau y}{m(x^2+y^2-1)^2}-\frac 1m\lambda^{\prime}\cos \theta+\frac 1m\lambda\theta^{\prime}\sin \theta,
                 \end{align*} together with  $\dot{x}\sin\theta=\dot{y}\cos\theta$ and boundary values as above.
                 
                 %$(x(0),y(0),\theta(0)),\qquad (x(T),y(T),\theta(T))$

%                 $(x^{\prime}(0),y^{\prime}(0),\theta^{\prime}(0)), \quad (x^{\prime}(T),y^{\prime}(T),\theta^{\prime}(T))$.}

\section{The variational obstacle avoidance problem on a Lie group}

Now we consider a Lie group $G$ endowed with a left-invariant  Riemannian metric $< \cdot , \cdot >$ defined by an inner product $\mathbb{I}$ on the Lie algebra $\mathfrak{g}$.
The Levi-Civita connection $\nabla$ induced by $< \cdot , \cdot >$ is an affine left-invariant connection and it is completely determined by its restriction to $\mathfrak{g}$ via left-translations. This restriction, denoted by $\stackrel{\mathfrak{g}}{\nabla}:\mathfrak{g}\times\mathfrak{g}\to\mathfrak{g}$,  is given by (see \cite{bookBullo} p. 271) $$\stackrel{\mathfrak{g}}{\nabla}_wu= \frac 12 [w,u]-\frac 12 \mathbb{I}^{\sharp}\left(\hbox{ad}_w^* \mathbb{I}^{\flat}(u)+\hbox{ad}_u^* \mathbb{I}^{\flat}(w)\right),$$ where \hbox{ad}$^{*}:\mathfrak{g}\times\mathfrak{g}^{*}\to\mathfrak{g}^{*}$ is the co-adjoint representation of $\mathfrak{g}$ on $\mathfrak{g}^{*}$ and where $\mathbb{I}^{\sharp}:\mathfrak{g}^{*}\to\mathfrak{g}$, $\mathbb{I}^{\flat}:\mathfrak{g}\to\mathfrak{g}^{*}$ are the associated isomorphisms to the inner product $\mathbb{I}$ (see \cite{Boothby} for instance).

If $u\in\mathfrak{g}$, its associated left-invariant vector field is given by $u_{L}(g)=T_{e}L_{g}(u)$ satisfying $u_{L}(e)=u\in T_{e}G$ where $L_g:G\to G$ denotes the left-translation map by $g$. If $u,v\in\mathfrak{g}$ it is possible to see that $\nabla_{w_L}u_L=(\stackrel{\mathfrak{g}}{\nabla}_wu)_L$ (see \cite{bookBullo} p. 273).

Let $x:I\subset\mathbb{R}\to G$ be a smooth curve on $G$. The \textit{body velocity} of $x$ is the curve $v:I\subset\mathbb{R}\to\mathfrak{g}$ defined by $\displaystyle{v(t)=T_{x(t)}L_{x(t)^{-1}}\left(\frac{dx}{dt}(t)\right)}$.

Let $\{e_1,\ldots,e_n\}$ be a basis of  $\mathfrak{g}$. Consider the body velocity of $x$ on the given basis, defined by $\displaystyle{v=\sum_{i=1}^n v_i e_i}$.
It follows that
\begin{equation}\label{admissibility}\frac{dx}{dt}(t)=T_eL_{x(t)}v(t)=\sum_{i=1}^n v_i(t) (e_i)_L(x(t)).\end{equation}

%For simplicity, omit the $t$ and  write $$\frac{dx}{dt}=T_eL_{x}v=\sum_{i=1}^n v_i (e_i)_L(x).$$
%Using this noation, we obtain
To write the equations determining necessary conditions for existence of a normal extremal, we must use the following formulas (see \cite{Altafini}, Section $7$ for more details)
%\begin{align*}
%&\frac{D^{2}x}{dt^{2}}=\nabla_{\frac{dx}{dt}}\frac{dx}{dt}=\nabla_{\frac{dx}{dt}}(\sum_{i=1}^n v_i (e_i)_L(x))=\sum_{i=1}^n v_i^{\prime} (e_i)_L(x)+\sum_{i=1}^n v_i \nabla_{\frac{dx}{dt}}((e_i)_L(x))=\\
%&\qquad=\sum_{i=1}^n v_i^{\prime} (e_i)_L(x)+\sum_{i,j=1}^n v_i v_j\nabla_{(e_j)_L(x)}((e_i)_L(x))=\sum_{i=1}^n v_i^{\prime} (e_i)_L(x)+\sum_{i,j=1}^n v_i v_j\Big(\nabla_{e_j}e_i\Big)_L(x)\\
%&\qquad=\Big(\sum_{i=1}^n v_i^{\prime} e_i+\sum_{i,j=1}^n v_i v_j\nabla_{e_j}e_i\Big)_L(x)=T_eL_{x}\Big(\sum_{i=1}^n v_i^{\prime} e_i+\sum_{i,j=1}^n v_i v_j\nabla_{e_j}e_i\Big).
%\end{align*}
%Denoting $\sum_{i=1}^n v_i^{\prime} e_i$ by $v^{\prime}$ and
\begin{align*}
&\stackrel{\mathfrak{g}}{\nabla}_vv=\sum_{i,j=1}^n v_i v_j\nabla_{e_j}e_i,\\% by $\stackrel{\mathfrak{g}}{\nabla}_vv$, %we can rewrite
&\frac{D^{2}x}{dt^{2}}=T_eL_{x}\Big(v^{\prime}+\stackrel{\mathfrak{g}}{\nabla}_vv\Big),\\
%In the same way, we conclude that
&\frac{D^{3}x}{dt^{3}}=T_eL_{x}\Big(v^{\prime \prime}+\stackrel{\mathfrak{g}}{\nabla}_{v^{\prime}}v+2 \stackrel{\mathfrak{g}}{\nabla}_vv^{\prime}+\stackrel{\mathfrak{g}}{\nabla}_v\stackrel{\mathfrak{g}}{\nabla}_vv\Big),\\
&\frac{D^{4}x}{dt^{4}}=T_eL_{x}\left(v'''+\stackrel{\mathfrak{g}}{\nabla}_{v''}v+3\stackrel{\mathfrak{g}}{\nabla}_{v^{\prime}}v^{\prime}+3 \stackrel{\mathfrak{g}}{\nabla}_vv''+\right.\\
&\qquad\qquad\qquad\left.\stackrel{\mathfrak{g}}{\nabla}_{v^{\prime}}\stackrel{\mathfrak{g}}{\nabla}_vv+2 \stackrel{\mathfrak{g}}{\nabla}_v\stackrel{\mathfrak{g}}{\nabla}_{v^{\prime}}v+3\stackrel{\mathfrak{g}}{\nabla}_v^2v^{\prime}+\stackrel{\mathfrak{g}}{\nabla}^{3}_vv\right),\\
&R\left(\frac{D^{2}x}{dt^{2}},\frac{dx}{dt}\right)\frac{dx}{dt}=T_eL_{x}\left(\mathfrak{R}(v^{\prime},v)v+\mathfrak{R}(\stackrel{\mathfrak{g}}{\nabla}_vv,v)v\right),\end{align*} where $\mathfrak{R}$ denotes the curvature tensor associated with $\stackrel{\mathfrak{g}}{\nabla}$.

 The equations giving rise to necessary conditions for the existence of normal extremals in the variational problem are, %(see Altafini for more details)
%\begin{equation*}
%\frac{D^{4}x}{dt^{4}}+R\left(\frac{D^2x}{dt^2},\frac{dx}{dt}\right)\frac{dx}{dt}+T_xL_{x^{-1}}(\mbox{grad }V_q(x))=\sigma
%\frac{D^2x}{dt^2}.
%\end{equation*}
\begin{align*}
0=&v'''+\stackrel{\mathfrak{g}}{\nabla}_{v''}v+3\stackrel{\mathfrak{g}}{\nabla}_{v^{\prime}}v^{\prime}+3 \stackrel{\mathfrak{g}}{\nabla}_vv''
+\stackrel{\mathfrak{g}}{\nabla}_{v^{\prime}}\stackrel{\mathfrak{g}}{\nabla}_vv\\
&+2\stackrel{\mathfrak{g}}{\nabla}_v\stackrel{\mathfrak{g}}{\nabla}_{v^{\prime}}v+
3\stackrel{\mathfrak{g}}{\nabla}_v^2v^{\prime}+\stackrel{\mathfrak{g}}{\nabla}_v^3v+\mathfrak{R}(v^{\prime},v)v-\sigma \stackrel{\mathfrak{g}}{\nabla}_vv\\&+\mathfrak{R}(\stackrel{\mathfrak{g}}{\nabla}_vv,v)v
-\sigma v^{\prime}+\frac{1}{2}T_xL_{x^{-1}}(\mbox{grad }V_q(x))\end{align*}
together with equation \eqref{admissibility}.

%Let $\{e_i\}_{i=\overline{1,n}}$ be a basis of  $\mathfrak{g}$. Consider the body velocity of $x$, $ v=\sum_{i=1}^n v_i (e_i)_L$ where $(e_i)_L$ are . We have  (see \cite{bookBullo} p. 276)
%
%$$\frac{dx}{dt}(t)=T_eL_{x(t)}v(t)=\sum_{i=1}^n v_i(t) (e_i)_L(x(t))$$
%For simplicity, omit the $t$ and  write $$\frac{dx}{dt}=T_eL_{x}v=\sum_{i=1}^n v_i (e_i)_L(x).$$
%In the same way, we have
%$$\frac{D^{2}x}{dt^{2}}=T_eL_{x}(v^{\prime}+\nabla_vv)$$
%where $\displaystyle{v^{\prime}=\sum_{i=1}^n v_i^{\prime} (e_i)_L(x)}$ and $\displaystyle{\nabla_vv= \sum_{i,j=1}^n v_i v_j\nabla_{(e_i)_L}(e_i)_L(x)}$,
%$$\frac{D^{3}x}{dt^{3}}=T_eL_{x}(v^{\prime \prime}+\nabla_{v^{\prime}}v+2 \nabla_vv^{\prime}+\nabla_v\nabla_vv)$$
%$$\frac{D^{4}x}{dt^{4}}=T_eL_{x}(v^{(3)}+\nabla_{v^{(2)}}v+3\nabla_{v^{\prime}}v^{\prime}+3 \nabla_vv^{(2)}+\nabla_{v^{\prime}}\nabla_vv+2 \nabla_v\nabla_{v^{\prime}}v+3\nabla_v^2v^{\prime}+\nabla_v^3v)$$
%
%$$R(\frac{D^{2}x}{dt^{2}},\frac{dx}{dt})\frac{dx}{dt}=...$$

%These formulas are important specially in $SE(3)$ case.

\subsection{ Example: planar rigid body on $SE(2)$.}

Consider $SE(2)$ endowed with a left-invariant metric representing the kinetic energy of a planar rigid body, defined by the inner product $\mathbb{I}=Je^1\otimes e^1+me^2\otimes e^2+me^3\otimes e^3$.
The Levi-Civita connection $\nabla$ induced by $< \cdot , \cdot >$ is left invariant and it is completely determined by its restriction to the Lie algebra $\mathfrak{se}(2)$, denoted by $\stackrel{\mathfrak{se(2)}}{\nabla}:\mathfrak{se}(2)\times\mathfrak{se}(2)\to\mathfrak{se}(2)$ and given by $$\stackrel{\mathfrak{se(2)}}{\nabla}_v w=-v_1(w_3e_2-w_2e_3)=\left(\begin{array}{c}
                                   0 \\
                                   -v_1w_3 \\
                                   v_1w_2
                                 \end{array}\right)
,$$ where $v=(v_1,v_2,v_3)$ and $w=(w_1,w_2,w_3)$ are the representative elements of $\mathfrak{se}(2)$ in $\mathbb{R}^3$  (see \cite{bookBullo} p. 279).
The curvature tensor is zero. Equations (\ref{admissibility}) become in
\begin{equation}\label{9} \theta^{\prime}=v_1, x^{\prime}=v_2\cos \theta - v_3 \sin \theta, y^{\prime}=v_2\sin \theta + v_3 \cos \theta\end{equation}

% As before we have
%\begin{equation}\label{bodycoord}
%\frac{dx}{dt}=T_eL_{x}v=\sum_{i=1}^3 v_i (e_i)_L(x).
%\end{equation}

%We can simply write $v=\left(\begin{array}{c}
   %                                v_1 \\
      %                             v_2 \\
         %                          v_3
            %                     \end{array}\right)$.

%$$\frac{D^{2}x}{dt^{2}}=T_eL_{x}(v^{\prime}+\nabla_vv)$$
%where $$ v^{\prime}+\nabla_vv=\left(\begin{array}{c}
%                                   v^{\prime}_1 \\
%                                   v^{\prime}_2-v_1v_3 \\
%                                   v^{\prime}_3+v_1v_2
%                                 \end{array}\right).$$
%
%$$\frac{D^{4}x}{dt^{4}}=T_eL_{x}(v^{(3)}+\nabla_{v^{(2)}}v+3\nabla_{v^{\prime}}v^{\prime}+3 \nabla_vv^{(2)}+\nabla_{v^{\prime}}\nabla_vv+2 \nabla_v\nabla_{v^{\prime}}v+3\nabla_v^2v^{\prime}+\nabla_v^3v)$$
%where
%$$v^{(3)}+\nabla_{v^{(2)}}v+3\nabla_{v^{\prime}}v^{\prime}+3 \nabla_vv^{(2)}+\nabla_{v^{\prime}}\nabla_vv+2 \nabla_v\nabla_{v^{\prime}}v+3\nabla_v^2v^{\prime}+\nabla_v^3v=$$
%$$\left(\begin{array}{c}
%                                  v_1^{(3)}\\
%                                   v_2^{(3)}-3 v_1{v_1^{\prime}}v_2-3v_1^2v_2^{\prime}+(v_1^3-v_1^{(2)})v_3-3v_1^{\prime}v_3^{\prime}-3v_1v_3^{(2)} \\
%                                   v_3^{(3)}-3 v_1{v_1^{\prime}}v_3-3v_1^2v_3^{\prime}+(v_1^3+v_1^{(2)})v_2+3v_1^{\prime}v_2^{\prime}+3v_1v_2^{(2)}
%                                 \end{array}\right)$$
%

The potential function $V: \mathrm{SE}(2)\to\mathbb{R}$ is given by
$$\displaystyle{V(g) = \frac{\tau}{(\|\hbox{Ad}_{g^{-1}}e_{1}\|^{2}-1)},}$$ where $
e_{1}$ is an element of the canonical basis for $\mathfrak{se}(2)$  and $\hbox{Ad}_{g}:\mathfrak{se}(2)\to\mathfrak{se}(2)$ is the adjoint representation of $SE(2)$ on $\mathfrak{se}(2)$ (see \cite{Gu}). A type of Euler-Poincar\'e equations can be obtained as in \cite{Gu}. We will study that approach in an ongoing work. Now, we restrict ourself to study the dynamics on $SE(2)\times\mathfrak{se}(2)\simeq\mathbb{R}^{2}\times S^{1}\times\mathbb{R}^{3}$. Using
$$
T_xL_{x^{-1}}(\mbox{grad } V)=-\frac {2V^2}{m\tau}(0,x\cos \theta+y\sin \theta,y\cos \theta-x\sin \theta)
$$ necessary conditions for existence of normal extremals of the variational obstacle avoidance problem are determined by
\begin{align*}
v_1'''&=\sigma v^{\prime}_1,\\
v_2'''&=3 v_1{v_1^{\prime}}v_2+3v_1^2v_2^{\prime}-(v_1^3-v_1'')v_3+3v_1^{\prime}v_3^{\prime}+3v_1v_3''\\&+\sigma (v^{\prime}_2-v_1v_3) +\frac {V^2}{m\tau} (x\cos \theta+y\sin \theta ),\\
v_3'''&=3 v_1{v_1^{\prime}}v_3+3v_1^2v_3^{\prime}+(v_1^3-v_1'')v_2-3v_1^{\prime}v_2^{\prime}\\&-3v_1v_2''+ \sigma(v^{\prime}_3+v_1v_2)  +\frac {V^2}{m\tau}(y\cos \theta-x\sin \theta),
\end{align*} together with equation \eqref{9}.

In the absence of obstacles, the equations reduce to the cubic polynomials in tension on $SE(2)$ \cite{SCC00}.
\begin{align*}
v_1'''&=\sigma v^{\prime}_1,\\
v_2'''&=3 v_1{v_1^{\prime}}v_2+3v_1^2v_2^{\prime}-(v_1^3-v_1'')v_3+3v_1^{\prime}v_3^{\prime}\\&+3v_1v_3''+\sigma (v^{\prime}_2+v_1v_3),\\
v_3'''&=3 v_1{v_1^{\prime}}v_3+3v_1^2v_3^{\prime}+(v_1^3-v_1'')v_2-3v_1^{\prime}v_2^{\prime}\\&-3v_1v_2''+ \sigma(v^{\prime}_3+v_1v_2),
\end{align*} together with equation \eqref{9}.
\section{Conclusions and future research}
We discussed obstacle avoidance variational problems on Riemannian manifolds and derived necessary conditions for the existence of normal extremals in the variational problem. Two different scenarios were studied: a general case on a Riemannian manifold and a sub-Riemannian problem. We also studied the variational obstacle avoidance problem on a Lie group.
	
The study of higher-order interpolation problems on arbitrary manifolds has attracted considerable interest and has been carried out systematically in the last decades by several authors. In a current work we incorporate interpolation points into the problem and we extend the results of this work and dynamic interpolation to variational obstacle avoidance problems.  In the last example, it is easy to verify that the potential function $V$ is $SO(2)$-invariant but not $SE(2)$-invariant and so, the potential function breaks the symmetry of the action functional \eqref{3.2} as in \cite{Gu}. This situation is also studied in our ongoing work.

Other interesting questions that we intend to study arrises in the situation when the Riemannian manifold is complete, and therefore we can connect points in the manifold by geodesics with the exponential function determining the geodesic distance. We will also consider extensions of dynamical interpolation to several obstacles into the picture of the problem as in \cite{Khatib1986}.

\section*{Acknowledgments}
The research of A. Bloch was supported by NSF grants DMS-1207693, DMS-1613819, INSPIRE-1343720 and the Simons Foundation. The research of M. Camarinha was partially supported by the Centre for Mathematics of the University of Coimbra -- UID/MAT/00324/2013, funded by the Portuguese
Government through FCT/MEC and co-funded by the European Regional
Development Fund through the Partnership Agreement PT2020. L. Colombo was supported by MINECO (Spain) grant MTM2016-76072-P.  L. C. wish to thank CMUC, Universidade de Coimbra for the hospitality received there where the main part of this work was developed.

\end{document}